\documentclass[12pt]{amsart}
\usepackage{amssymb, hyperref}

\newcommand{\Bgp}{{\Z^\N}}

\long\def\forget#1\forgotten{}
\newcommand{\issuenumber}{31}
\newcommand{\issuemonth}{April}
\newcommand{\issueyear}{2011}

\newcommand{\ed}{
\newpage

\section{Unsolved problems from earlier issues}

\begin{issue}
Is $\binom{\Omega}{\Gamma}=\binom{\Omega}{\Tau}$?
\end{issue}

\begin{issue}
Is $\ufin(\cO,\Omega)=\sfin(\Gamma,\Omega)$?
And if not, does $\ufin(\cO,\Gamma)$ imply
$\sfin(\Gamma,\Omega)$?
\end{issue}

\stepcounter{issue}

\begin{issue}
Does $\sone(\Omega,\Tau)$ imply $\ufin(\Gamma,\Gamma)$?
\end{issue}

\begin{issue}
Is $\fp=\fp^*$? (See the definition of $\fp^*$ in that issue.)
\end{issue}

\begin{issue}
Does there exist (in ZFC) an uncountable set satisfying $\sfin(\cB,\cB)$?
\end{issue}

\stepcounter{issue}

\begin{issue}
Does $X \nin \NON(\cM)$ and $Y\nin\mathsf{D}$ imply that
$X\cup Y\nin \COF(\cM)$?
\end{issue}

\begin{issue}[CH]
Is $\split(\Lambda,\Lambda)$ preserved under finite unions?
\end{issue}

\begin{issue}
Is $\cov(\cM)=\fo$? (See the definition of $\fo$ in that issue.)
\end{issue}

\stepcounter{issue}

\begin{issue}
Could there be a Baire metric space $M$ of weight $\aleph_1$ and a partition
$\mathcal{U}$ of $M$ into $\aleph_1$ meager sets where for each ${\mathcal U}'\subset\mathcal U$,
$\bigcup {\mathcal U}'$ has the Baire property in $M$?
\end{issue}

\stepcounter{issue} 

\begin{issue}
Does there exist (in ZFC) a set of reals $X$ of cardinality $\fd$ such that all
finite powers of $X$ have Menger's property $\sfin(\cO,\cO)$?
\end{issue}

\begin{issue}
Can a Borel non-$\sigma$-compact group be generated by a Hurewicz subspace?
\end{issue}

\begin{issue}[MA]
Is there $X\sbst\bbR$ of cardinality continuum, satisfying $\sone(\BO,\BG)$?
\end{issue}

\begin{issue}[CH]
Is there a totally imperfect $X$ satisfying $\ufin(\cO,\Gamma)$
that can be mapped continuously onto $\Cantor$?
\end{issue}

\begin{issue}[CH]
Is there a Hurewicz $X$ such that $X^2$ is Menger but not Hurewicz?
\end{issue}

\begin{issue}
Does the Pytkeev property of $C_p(X)$ imply that $X$ has Menger's property?
\end{issue}

\begin{issue}
Does every hereditarily Hurewicz space satisfy $\sone(\BG,\BG)$?
\end{issue}

\begin{issue}[CH]
Is there a Rothberger-bounded $G\le\Bgp$ such that $G^2$ is not Menger-bounded?
\end{issue}

\begin{issue}
Let $\cW$ be the van der Waerden ideal.
Are $\cW$-ultrafilters closed under products?
\end{issue}

\begin{issue}
Is the $\delta$-property equivalent to the $\gamma$-property $\binom{\Omega}{\Gamma}$?
\end{issue}

\stepcounter{issue}

\stepcounter{issue}

\general\end{document}}

\newcommand{\Cantor}{{\{0,1\}^\N}}

\newcommand{\fb}{\mathfrak{b}}

\newcommand{\fc}{\mathfrak{c}}
\newcommand{\fd}{\mathfrak{d}}
\newcommand{\fp}{\mathfrak{p}}
\newcommand{\fs}{\mathfrak{s}}

\newcommand{\NON}{{\mathsf   {NON}}}
\newcommand{\COF}{{\mathsf   {COF}}}

\newcommand{\cN}{\mathcal{N}}

\newcommand{\cM}{\mathcal{M}}

\newcommand{\cov}{\mathsf{cov}}

\newcommand{\bbR}{\mathbb{R}}

\newcommand{\EdNote}[1]{\par\medskip\noindent\textbf{#1.}}
\newcommand{\fo}{\mathfrak{od}}

\newcommand{\w}{\omega}

\renewcommand{\split}{\mathsf{Split}}
\newcommand{\bq}{\begin{quote}}
\newcommand{\eq}{\end{quote}}
\newcommand{\cO}{\mathcal{O}}
\newcommand{\cK}{\mathcal{K}}
\newcommand{\cB}{\mathcal{B}}
\newcommand{\BG}{\cB_\Gamma}

\newcommand{\BO}{\cB_\Omega}

\newcommand{\sone}{\mathsf{S}_1}    \newcommand{\sfin}{\mathsf{S}_\mathrm{fin}}

\newcommand{\ufin}{\mathsf{U}_\mathrm{fin}}

\newcommand{\nin}{\not\in}
\newcommand{\cF}{\mathcal{F}}

\newcommand{\cU}{\mathcal{U}}

\newcommand{\cW}{\mathcal{W}}

\newcommand{\N}{\mathbb{N}}

\newcommand{\Z}{\mathbb{Z}}

\newcommand{\sbst}{\subseteq}
\newcommand{\by}[2]{\par\hfill\emph{#1}, #2}
\newcommand{\nby}[1]{\par\hfill\emph{#1}}
\newcommand{\Tau}{\mathrm{T}}
\newcommand{\CE}{\textsc{CE}}

\newtheorem{thm}{Theorem}[section]
\newcommand{\bthm}{\begin{thm}} \newcommand{\ethm}{\end{thm}}
\newtheorem{prop}[thm]{Proposition}
\newcommand{\bprp}{\begin{prop}} \newcommand{\eprp}{\end{prop}}
\newtheorem{fact}[thm]{Fact}
\newcommand{\bfct}{\begin{fact}} \newcommand{\efct}{\end{fact}}
\newtheorem{prob}[thm]{Problem}
\newcommand{\bprb}{\begin{prob}} \newcommand{\eprb}{\end{prob}}
\newtheorem{lem}[thm]{Lemma}
\newcommand{\blem}{\begin{lem}} \newcommand{\elem}{\end{lem}}
\newtheorem{claim}[thm]{Claim}
\newcommand{\bclm}{\begin{claim}} \newcommand{\eclm}{\end{claim}}
\newtheorem{cor}[thm]{Corollary}
\newcommand{\bcor}{\begin{cor}} \newcommand{\ecor}{\end{cor}}
\newtheorem{conj}[thm]{Conjecture}
\newcommand{\bcnj}{\begin{conj}} \newcommand{\ecnj}{\end{conj}}
\theoremstyle{definition}
\newtheorem{defn}[thm]{Definition}
\newcommand{\bdfn}{\begin{defn}} \newcommand{\edfn}{\end{defn}}
\theoremstyle{remark}
\newtheorem{rem}[thm]{Remark}
\newcommand{\brem}{\begin{rem}} \newcommand{\erem}{\end{rem}}
\newtheorem{cnv}[thm]{Convention}
\newcommand{\bcnv}{\begin{cnv}} \newcommand{\ecnv}{\end{cnv}}
\newtheorem{exam}[thm]{Example}
\newcommand{\bexm}{\begin{exam}} \newcommand{\eexm}{\end{exam}}
\newtheorem{issue}{Issue}

\newcommand{\bpf}{\begin{proof}} \newcommand{\epf}{\end{proof}}
\newcommand{\be}{\begin{enumerate}}
\newcommand{\ee}{\end{enumerate}}
\newcommand{\bi}{\begin{itemize}}
\newcommand{\ei}{\end{itemize}}
\newcommand{\itm}{\item}

\newcommand{\general}{\small\vfill\par\noindent\hrulefill\par
\noindent\textbf{Previous issues.} The previous issues of this
bulletin are available online at\\
\url{http://front.math.ucdavis.edu/search?\&t=\%22SPM+Bulletin\%22}
\\[0.1cm]
\textbf{Contributions.} Announcements, discussions, and open problems should be emailed
to \texttt{tsaban@math.biu.ac.il}\\[0.1cm]
\textbf{Subscription.}
To receive this bulletin (free) to your e-mailbox, e-mail us.
}

\newcommand{\arXivl}[4]{\subsection{#2}{#4}\par\hfill{\arx{#1}}\par\hfill\emph{#3}}
\newcommand{\arXiv}[3]{\subsection{#2}\mbox{}\par\hfill{\arx{#1}}\par\hfill\emph{#3}}

\newcommand{\arx}[1]{\url{http://arxiv.org/abs/#1}}

\title[$\mathcal{SPM}$ Bulletin \textbf{\issuenumber} (\issuemonth{} \issueyear)]{%
$\mathcal{SPM}$ Bulletin\\[0.5cm]
Issue number \issuenumber: \issuemonth{} \issueyear{} \CE{}}

\begin{document}
\maketitle


\section{Editor's note}

A recent series of papers of Franklin Tall on selective properties (SPM),
some of which announced below, is noteworthy.

Greetings to Vladimir Tkachuk for the publication of his new book, announced below.

\medskip

\by{Boaz Tsaban}{tsaban@math.biu.ac.il}

\hfill \texttt{http://www.cs.biu.ac.il/\~{}tsaban}

\section{A new book on $C_p$-theory}
Dear Colleagues: This message is  to  inform you that  my  book entitled \textbf{A $C_p$-theory Problem Book}
has already been  published  in  Springer.  At  the Springer's page\\
\url{http://www.springer.com/mathematics/geometry/book/978-1-4419-7441-9}\\
you  can see its contents, sample pages and preface. The book has
500  problems   with   complete   solutions   and  constitutes  a
self-contained introduction to $C_p$-theory and  General  Topology.
However,  it  also  contains  research topics and deep results in
general topology and $C_p$-theory. To  mention  just a few, in this
book you can find:
\be
\itm Ten properties equivalent to paracompactness.
\itm The Stone-Weierstrass theorem for compact spaces.
\itm Theorems on cardinal functions in linearly ordered spaces.
\itm Introduction to the theory of realcompact spaces and Dieudonne
complete spaces.
\itm Shapirovsky's deep  theorem  that states  that every compact
space of countable tightness has a point-countable $\pi$-base.
\itm The theorem that every continuous map on a product  of  second
countable spaces depends on countably many coordinates.
\itm Arhagel'skii's  theorem  on  cardinality  of  first countable
compact spaces.
\itm Arhagel'skii's theorem on  tightness  and free sequences in
compact spaces.
\itm The  theorem  states  that  every  dyadic  compact  space  of
countable tightness is metrizable.
\itm Shakhmatov's example of an infinite space X such that $C_p(X)$
is $\sigma$-pseudocom\-pact.
\ee
and many other concepts, facts an theorems together with 100 open
problems  in  $C_p$-theory  and  a bibliography of 200 items.

With best regards,

\nby{Vladimir Tkachuk}

\section{Long announcements}

\arXivl{1007.4309}
{Elementary submodels in infinite combinatorics}
{Lajos Soukup}
{We show that usage of elementary submodels is a simple but powerful method to
prove theorems, or to simplify proofs in infinite combinatorics. First we
introduce all the necessary concepts of logic, then we prove classical theorems
using elementary submodels. We also present a new proof of Nash-Williams's
theorem on cycle-decomposition of graphs, and finally we obtain some new
decomposition theorems by eliminating GCH from some proofs concerning
bond-faithful decompositions of graphs.}

\arXivl{1007.5368}
{Monotone hulls for $\cN\cap\cM$}
{Andrzej Roslanowski and Saharon Shelah}
{Using the method of decisive creatures (see Kellner and Shelah [KrSh:872])
we show the consistency of ``there is no increasing
 $\omega_2$--chain of Borel sets and ${\rm non}({\mathcal N})= {\rm
   non}({\mathcal M})=\omega_2=2^\omega$''. Hence, consistently, there are
 no monotone Borel hulls for the ideal ${\mathcal M}\cap {\mathcal
   N}$. This answers questions of Balcerzak and Filipczak.
 Next we use FS iteration with partial memory to show that
 there may be monotone Borel hulls for the ideals $\cM,\cN$ even if they
 are not generated by towers.
}

\EdNote{Editor's note}
Steprans and I observed long ago (answering a question A. Krawczyk asked us on our way
back from the first European Set Theory conference)
that in the Cohen model, there is no increasing
$\aleph_2$-chain of Borel sets. This can be proved in the same way as Kunen's classical (thesis) argument
that there is no tower of length $\aleph_2$ in the Cohen model. This simple observation was not published.
The above result is much stronger, of course.

\arXivl{1008.4739}
{Continuous maps on Aronszajn trees}
{Kenneth Kunen, Jean Larson, and Juris Stepr\=ans}
{Assuming Jenson's principle diamond: Whenever $B$ is a totally imperfect set of
real numbers, there is special Aronszajn tree with no continuous order
preserving map into $B$.}

\arXivl{1009.0065}
{The Filter Dichotomy and medial limits}
{Paul B. Larson}
{The \emph{Filter Dichotomy} says that every uniform nonmeager filter on the
integers is mapped by a finite-to-one function to an ultrafilter. The
consistency of this principle was proved by Blass and Laflamme. A function
between topological spaces is \emph{universally measurable} if the preimage of
every open subset of the codomain is measured by every Borel measure on the
domain. A \emph{medial limit} is a universally measurable function from
$\mathcal{P}(\omega)$ to the unit interval $[0,1]$ which is finitely additive for
disjoint sets, and maps singletons to $0$ and $\omega$ to 1. Christensen and
Mokobodzki independently showed that the Continuum Hypothesis implies the
existence of medial limits. We show that the Filter Dichotomy implies that
there are no medial limits.}

\arXivl{1009.0818}
{Bernstein sets and $\kappa$-coverings}
{J. Kraszewski, R. Ralowski, P. Szczepaniak and S. Zeberski}
{In this paper we study a notion of a $\kappa$-covering in connection with
Bernstein sets and other types of nonmeasurability. Our results correspond to
those obtained by Muthuvel and Nowik. We consider also other types of
coverings.}

\arXivl{1009.3683}
{Ideal games and Ramsey sets}
{Carlos Di Prisco, Jose G. Mijares, Carlos Uzcategui}
{It is shown that Matet's characterization of $\mathcal{H}$-Ramseyness
relative to a selective coideal $\mathcal{H}$, in terms of games of Kastanas,
still holds if we consider semiselectivity instead of selectivity. Moreover, we
prove that a coideal $\mathcal{H}$ is semiselective if and only if Matet's
game-theoretic characterization of $\mathcal{H}$-Ramseyness holds. This gives a
game-theoretic counter part to a theorem of Farah, asserting that a coideal
$\mathcal{H}$ is semiselective if and only if the family of
$\mathcal{H}$-Ramsey subsets of $\mathbb{N}^{[\infty]}$ coincides with the
family of those sets having the $Exp(\mathcal{H})$-Baire property. Finally, we
show that under suitable assumptions, semiselectivity is equivalent to the
Fr\'echet-Urysohn property.}

\arXivl{1010.0999}
{Nonmeasurable unions of sets and continuity of group representations}
{Julia Kuznetsova}
{Let $G$ be a locally compact group, and let $U$ be its unitary representation on
a Hilbert space $H$. Endow the space $L(H)$ of linear bounded operators on $H$ with
weak operator topology. We prove that if $U$ is a measurable map from $G$ to $L(H)$
then it is continuous. This result was known before for separable $H$. To prove
this, we generalize a known theorem on nonmeasuralbe unions of point finite
families of null sets. We prove also that the following statement is consistent
with ZFC: every measurable homomorphism from a locally compact group into any
topological group is continuous. This relies, in turn, on the following
theorem: it is consistent with ZFC that for every null set $S$ in a locally
compact group there is a set $A$ such that $AS$ is non-measurable.}

\arXivl{1010.1226}
{On weakly tight families}
{Dilip Raghavan and Juris Stepr\=ans}
{Using ideas from Shelah's recent proof that a completely separable maximal
almost disjoint family exists when $\fc < {\aleph}_{\omega}$, we construct a
weakly tight family under the hypothesis $\fs \leq \fb < {\aleph}_{\omega}$. The
case when $\fs < \fb$ is handled in ZFC and does not require $\fb <
{\aleph}_{\omega}$, while an additional PCF type hypothesis, which holds when
$\fb < {\aleph}_{\omega}$ is used to treat the case $\fs = \fb$. The notion of a
weakly tight family is a natural weakening of the well studied notion of a
Cohen indestructible maximal almost disjoint family. It was introduced by
Hru{\v{s}}{\'a}k and Garc{\'{\i}}a Ferreira, who applied it to the
Kat\'etov order on almost disjoint families.}

\arXivl{1010.1359}
{Partitions of groups and matroids into independent subsets}
{Taras Banakh, Igor Protasov}
{Can the real line with removed zero be covered by countably many linearly
(algebraically) independent subsets over the field of rationals? We use a
matroid approach to show that an answer is ``Yes'' under the Continuum
Hypothesis, and ``No'' under its negation.}

\arXivl{1010.2474}
{On $M$-separability of countable spaces and function spaces}
{Du\v{s}an Repov\v{s} and Lyubomyr Zdomskyy}
{We study $M$-separability as well as some other combinatorial versions of
separability. In particular, we show that the set-theoretic hypothesis $\fb=\fd$
implies that the class of selectively separable spaces is not closed under
finite products, even for the spaces of continuous functions with the topology
of pointwise convergence. We also show that there exists no maximal $M$-separable
countable space in the model of Frankiewicz, Shelah, and Zbierski in which all
closed $P$-subspaces of $\omega^*$ admit an uncountable family of nonempty open mutually
disjoint subsets. This answers several questions of Bella, Bonanzinga, Matveev,
and Tkachuk.}

\arXivl{1010.3368}
{Bornologies, selection principles and function spaces}
{Agata Caserta, Giuseppe Di Maio and Ljubisa D.R. Kocinac}
{We study some closure-type properties of function spaces endowed with the new
topology of strong uniform convergence on a bornology introduced by Beer and
Levy in 2009. The study of these function spaces was initiated elsewhere.
The properties we study are related to selection principles.}

\arXivl{1011.1869}
{Rothberger bounded groups and Ramsey theory}
{Marion Scheepers}
{We show that: 1. Rothberger bounded subgroups of sigma-compact groups are
characterized by Ramseyan partition relations. 2. For each uncountable cardinal
$\kappa$ there is a ${\sf T}_0$ topological group of cardinality $\kappa$ such
that ONE has a winning strategy in the point-open game on the group and the
group is not a subspace of any sigma-compact space. 3. For each uncountable
cardinal $\kappa$ there is a ${\sf T}_0$ topological group of cardinality
$\kappa$ such that ONE has a winning strategy in the point-open game on the
group and the group is $\sigma$-compact.}

\arXivl{1011.1031}
{On the length of chains of proper subgroups covering a topological group}
{Taras Banakh, Du\v{s}an Repov\v{s}, Lyubomyr Zdomskyy}
{We prove that if an ultrafilter $L$ is
not coherent to a $Q$-point, then
 each  analytic non-$\sigma$-bounded
 topological group $G$ admits an increasing
chain $\langle G_\alpha:\alpha<\fb(L)\rangle$ of its proper
  subgroups such that:
\be
\itm $\bigcup_{\alpha}G_\alpha=G$; and
\itm For every $\sigma$-bounded subgroup $H$ of $G$ there exists $\alpha$ such
 that $H\subset G_\alpha$.
\ee
In case of the group of all permutations of $\w$
with the topology inherited from $\w^\w$, this
improves upon earlier results of  S.~Thomas.
}

\arXivl{1011.3574}
{Core compactness and diagonality in spaces of open sets}
{Francis Jordan and Frederic Mynard}
{We investigate when the space $\mathcal O_X$ of open subsets of a topological
space $X$ endowed with the Scott topology is core compact. Such conditions turn
out to be related to infraconsonance of $X$, which in turn is characterized in
terms of coincidence of the Scott topology of $\mathcal O_X\times\mathcal O_X$
with the product of the Scott topologies of $\mathcal O_X$ at $(X,X)$. On the
other hand, we characterize diagonality of $\mathcal O_X$ endowed with the
Scott convergence and show that this space can be diagonal without being
pretopological. New examples are provided to clarify the relationship between
pretopologicity, topologicity and diagonality of this important convergence
space.}

\arXivl{1011.3586}
{Very I-favorable spaces}
{A. Kucharski, Sz. Plewik and V. Valov}
{We prove that a Hausdorff space $X$ is very $\mathrm I$-favorable if and only
if $X$ is the almost limit space of a $\sigma$-complete inverse system
consisting of (not necessarily Hausdorff) second countable spaces and
surjective d-open bonding maps. It is also shown that the class of Tychonoff
very $\mathrm I$-favorable spaces with respect to the co-zero sets coincides
with the d-openly generated spaces.}

\arXivl{1011.4555}
{Topological classification of zero-dimensional $M_\omega$-groups}
{Taras Banakh}
{A topological group $G$ is called an $M_\omega$-group if it admits a
countable cover $\cK$ by closed metrizable subspaces of $G$ such that a subset
$U$ of $G$ is open in $G$ if and only if $U\cap K$ is open in $K$ for every
$K\in\cK$. It is shown that any two non-metrizable uncountable separable
zero-dimenisional $M_\omega$-groups are homeomorphic. Together with Zelenyuk's
classification of countable $k_\omega$-groups this implies that the topology of
a non-metrizable zero-dimensional $M_\omega$-group $G$ is completely determined
by its density and the compact scatteredness rank $r(G)$ which, by definition,
is equal to the least upper bound of scatteredness indices of scattered compact
subspaces of $G$.}

\arXivl{1012.1094}
{On the Menger covering property and $D$-spaces}
{Du\v{s}an Repov\v{s} and Lyubomyr Zdomskyy}
{The main results of this note are: It is consistent that every subparacompact
space $X$ of size $\omega_1$ is a $D$-space; If there exists a Michael space,
then all productively Lindel\"of spaces have the Menger property, and,
therefore, are $D$-spaces; and
 Every locally $D$-space which admits a $\sigma$-locally finite cover by
Lindel\"of spaces is a $D$-space.}

\arXivl{1012.2522}
{On meager function spaces, network character and meager convergence in topological spaces}
{Taras Banakh, Volodymyr Mykhaylyuk, Lyubomyr Zdomskyy}
{For a non-isolated point $x$ of a topological space $X$ the network character
$nw_\chi(x)$ is the smallest cardinality of a family of infinite subsets of $X$
such that each neighborhood $O(x)$ of $x$ contains a set from the family. We
prove that (1) each paracompact space $X$ admitting a closed map onto a
non-discrete Frechet-Urysohn space contains a non-isolated point $x$ with
countable network character; (2) for each point $x\in X$ with countable
character there is an injective sequence in $X$ that $\cF$-converges to $x$ for
some meager filter $\cF$ on $\omega$; (3) if a functionally Hausdorff space $X$
contains an $\cF$-convergent injective sequence for some meager filter $\cF$,
then for every $T_1$-space $Y$ that contains two non-empty open sets with
disjoint closures, the function space $C_p(X,Y)$ is meager.}

\arXivl{1012.3422}
{Independently axiomatizable $L_{\omega_1,\omega}$ theories}
{Greg Hjorth, Ioannis Souldatos}
{In partial answer to a question posed by Arnie Miller
(\texttt{http://www.math.wisc.edu/\~{}miller/res/problem.pdf}) and X.
Caicedo, we obtain
sufficient conditions for an $L_{\omega_1,\omega}$ theory to have an independent
axiomatization. As a consequence we obtain two corollaries: The first, assuming
Vaught's Conjecture, every $L_{\omega_1,\omega}$ theory in a countable
language has
an independent axiomatization.

The second, this time outright in ZFC, every
intersection of a family of Borel sets can be formed as the intersection of a
family of independent Borel sets.

J. Symbolic Logic \textbf{74} (2009),  1273-1286.
}

\arXivl{1012.3966}
{Order-theoretic properties of bases in topological spaces, I}
{Menachem Kojman, David Milovich and Santi Spadaro}
{We study some cardinal invariants of an order-theoretic fashion on products
and box products of topological spaces. In particular, we concentrate on the
Noetherian type (Nt), defined by Peregudov in the 1990s. Some highlights of our
results include: 1) There are spaces $X$ and $Y$ such that $Nt(X \times Y) <
\min\{Nt(X), Nt(Y)\}$. 2) In several classes of compact spaces, the Noetherian
type is preserved by their square and their dense subspaces. 3) The Noetherian
type of some countably supported box products cannot be determined in ZFC. In
particular, it is sensitive to square principles and some Chang Conjecture
variants. 4) PCF theory can be used to provide ZFC upper bounds to Noetherian
type on countably supported box products. The underlying combinatorial notion
is a weakening of Shelah's freeness.}

\arXivl{1012.4338}
{Quasi-selective and weakly Ramsey ultrafilters}
{Marco Forti}
{Selective ultrafilters are characterized by many equivalent properties, in
particular the Ramsey property that every finite coloring of unordered pairs
of integers has a homogeneous set in $\cU$, and the equivalent property that every
function is nondecreasing on some set in $\cU$. Natural weakenings of these
properties led to the inequivalent notions of weakly Ramsey and of
quasi-selective ultrafilter, introduced and studied in earlier works.
$\cU$ is weakly Ramsey if for every finite coloring of unordered
pairs of integers there is a set in $\cU$ whose pairs share only two colors, while
$\cU$ is $f$-quasi-selective if every function $g < f$ is nondecreasing on some set in
$\cU$.
In this paper we consider the relations between various natural cuts of the
ultrapowers of $\N$ modulo weakly Ramsey and $f$-quasi-selective ultrafilters. In
particular we characterize those weakly Ramsey ultrafilters that are isomorphic
to a quasi-selective ultrafilter.}

\arXivl{1101.2754}
{Topologies on groups determined by sets of convergent sequences}
{S. Gabriyelyan}
{A Hausdorff topological group $(G,\tau)$ is called a $s$-group and $\tau$ is
called a $s$-topology if there is a set $S$ of sequences in $G$ such that
$\tau$ is the finest Hausdorff group topology on $G$ in which every sequence of
$S$ converges to the unit. The class $\mathbf{S}$ of all $s$-groups contains
all sequential Hausdorff groups and it is finitely multiplicative. A quotient
group of a $s$-group is a $s$-group. For non-discrete (Abelian) topological
group $(G,\tau)$ the following three assertions are equivalent: 1) $(G,\tau)$
is a $s$-group, 2) $(G,\tau)$ is a quotient group of a Graev-free (Abelian)
topological group over a Fr\'{e}chet-Urysohn Tychonoff space, 3) $(G,\tau)$ is
a quotient group of a Graev-free (Abelian) topological group over a sequential
Tychonoff space.}

\arXivl{1101.4615}
{Variations of selective separability II: discrete sets and the influence of convergence and maximality}
{Angelo Bella, Mikhail Matveev, Santi Spadaro}
{A space $X$ is called selectively separable(R-separable) if for every sequence of dense subspaces $(D_n : n\in\omega)$ one can pick finite (respectively, one-point) subsets $F_n\subset D_n$ such that $\bigcup_{n\in\omega}F_n$ is dense in $X$. These properties are much stronger than separability, but are equivalent to it in the presence of certain convergence properties. For example, we show that every Hausdorff separable radial space is R-separable and note that neither separable sequential nor separable Whyburn spaces have to be selectively separable. A space is called \emph{d-separable} if it has a dense $\sigma$-discrete subspace. We call a space $X$ D-separable if for every sequence of dense subspaces $(D_n : n\in\omega)$ one can pick discrete subsets $F_n\subset D_n$ such that $\bigcup_{n\in\omega}F_n$ is dense in $X$. Although $d$-separable spaces are often also $D$-separable (this is the case, for example, with linearly ordered $d$-separable or stratifiable spaces), we offer three examples of countable non-$D$-separable spaces. It is known that d-separability is preserved by arbitrary products, and that for every $X$, the power $X^{d(X)}$ is d-separable. We show that D-separability is not preserved even by finite products, and that for every infinite $X$, the power $X^{2^{d(X)}}$ is not D-separable. However, for every $X$ there is a $Y$ such that $X\times Y$ is D-separable. Finally, we discuss selective and D-separability in the presence of maximality. For example, we show that (assuming ${\mathfrak d}=\mathfrak c$) there exists a maximal regular countable selectively separable space, and that (in ZFC) every maximal countable space is D-separable (while some of those are not selectively separable). However, no maximal space satisfies the natural game-theoretic strengthening of D-separability.}

\arXivl{1103.5716}
{Star-covering properties: generalized $\Psi$-spaces, countability conditions, reflection}
{L. P. Aiken}
{We investigate star-covering properties of $\Psi$-like spaces. We show star-Lindel\"ofness is reflected by open perfect mappings. In addition, we offer a new equivalence of CH.}

\arXivl{1104.1759}
{On productively Lindel\"of spaces}
{Franklin D. Tall, Boaz Tsaban}
{The class of spaces such that their product with every Lindel\"of space is Lindel\"of is not well-understood. We prove a number of new results concerning such productively Lindel\"of spaces with some extra property, mainly assuming the Continuum Hypothesis.}

\arXivl{1104.2794}
{Productively Lindelof spaces may all be D}
{Franklin D. Tall}
{We give easy proofs that a) the Continuum Hypothesis implies that if the product of $X$ with every Lindel\"of space is Lindel\"of, then $X$ is a D-space, and b) Borel's Conjecture implies every Rothberger space is Hurewicz.}

\EdNote{Remark}
Reader's puzzled by (b) of the last announcement should notice that it deals with \emph{arbitrary},
not necessarily metrizable, topological spaces.

\arXiv{1104.2793}
{Lindel\"of spaces which are indestructible, productive, or D}
{Franklin D. Tall, Leandro F. Aurichi}
{We discuss relationships in Lindel\"of spaces among the properties ``indestructible'', ``productive'', ``D'', and related properties.}

\arXivl{1104.2796}
{Set-theoretic problems concerning Lindel\"of spaces}
{Franklin D. Tall}
{I survey problems concerning Lindel\"of spaces which have partial set-
theoretic solutions.}

\section{Short announcements}\label{RA}

\arXivl{1007.2266}
{The strong splitting number}
{Shimon Garti and Saharon Shelah}

\arXiv{1006.3816}
{On strongly summable ultrafilters}
{Peter Krautzberger}

\arXiv{1007.1666}
{Guessing clubs for aD, non D-spaces}
{Daniel Soukup}

\arXiv{1007.3367}
{Polish topometric groups}
{Ita\"i Ben Yaacov, Alexander Berenstein, Julien Melleray}

\arXiv{1007.2693}
{Large weight does not yield an irreducible base}
{Saharon Shelah}

\arXiv{1007.4034}
{All automorphisms of all Calkin algebras}
{Ilijas Farah}

\arXiv{1010.0327}
{Forcing $\square_{\omega_1}$ with finite conditions}
{Gregor K. Dolinar and Mirna D\v{z}amonja}

\arXiv{1010.3329}
{A decomposition theorem for compact groups with application to supercompactness}
{Wies{\l}aw Kubi\'s, S{\l}awomir Turek}

\arXiv{1011.1524}
{Productivity of sequences with respect to a given weight function}
{Dikran Dikranjan, Dmitri Shakhmatov, Jan Sp\v{e}v\'ak}

\arXiv{1011.3530}
{Dual topologies on non-abelian groups}
{Mar\'ia V. Ferrer, Salvador Hern\'andez}

\arXiv{1011.2089}
{Quasi-selective ultrafilters and asymptotic numerosities}
{Mauro Di Nasso and Marco Forti}

\arXiv{1011.4554}
{Topologies on groups determined by sequences: Answers to several questions of I.Protasov and E.Zelenyuk}
{Taras Banakh}

\arXiv{1012.0596}
{A direct proof of the five element basis theorem}
{Boban Velickovic, Giorgio Venturi}

\arXiv{1012.2040}
{The combinatorial essence of supercompactness}
{Christoph Wei{\ss}}

\arXiv{1012.2046}
{On the consistency strength of the proper forcing axiom}
{Matteo Viale, Christoph Wei{\ss}}

\arXiv{1012.4177}
{A Kronecker-Weyl theorem for subsets of abelian groups}
{Dmitri Shakhmatov, Dikran Dikranjan}

\arXiv{1012.4954}
{Finite basis for analytic strong $n$-gaps}
{Antonio Avil\'es, Stevo Todorcevic}

\arXiv{1101.2756}
{Pontryagin duality for Abelian $s$- and $sb$-groups}
{S. Gabriyelyan}

\arXiv{1101.4504}
{Pontryagin duality in the class of precompact Abelian groups and the Baire property}
{Montserrat Bruguera and Mikhail Tkachenko}

\arXiv{1102.5077}
{Metrization criteria for compact groups in terms of their dense subgroups}
{Dikran Dikranjan, Dmitri Shakhmatov}

\arXiv{1104.2795}
{PFA(S)[S] and the Arhangel'skii-Tall problem}
{Franklin D. Tall}

\ed